\documentclass[12pt]{amsart}

\usepackage{amsmath}

\font\smallcaps=cmcsc10
\def\maple{{\smallcaps maple\ }} 
\def\maplep{{\smallcaps maple.\ }}

\def\carrot{${}^\wedge$}

\def\ssn{\vskip 3pt\noindent}
\def\sssn{\vskip 0pt\noindent}
\def\sn{\vskip 10pt\noindent}

\def\qed{\penalty-250%
\hbox to 6pt{\hfill}\hfill\llap{%
\vbox{\hrule\hbox{\vrule height6pt\hskip6pt\vrule}\hrule}}}
\def\t#1{{\tt #1}}
\def\EX#1#2{\sn
\advance\leftskip by 20pt \advance\rightskip by 20pt 
\sn
{{\bf Exercise #1}.} \quad #2
\sn
\advance\leftskip by -20pt \advance\rightskip by -20pt 
\sn
}
\def\sec#1{(\ref{sec-#1})}
\def\eqn#1{(\ref{eq:#1})}
\def\qinf#1#2{(q^{#1};q^{#2})_\infty}
\def\receivedline{\relax}
\def\dedication#1{\receivedline\vskip4pt
{\normalsize\begin{center}#1\end{center}\vskip1sp}}

\begin{document}

\title[Qseries maple tutorial]  %
      {A $q$-product tutorial\\
	for a\\
       a $q$-series maple package}

\date{Tuesday, December 15, 1998}

\author{Frank Garvan}
\address{Department of Mathematics \\
         University of Florida  \\
         Gainesville, Florida 32611} 
\email{frank@math.ufl.edu}
\thanks{This research was supported by 
   the NSF under grant number~ DMS-9870052.} 


\begin{abstract}
	This is a tutorial for using a new $q$-series 
	Maple package. The package includes facilities for
	conversion between $q$-series and $q$-products
	and finding algebraic relations between $q$-series.
	Andrews found an algorithm for converting a $q$-series
	into a product. We provide an implementation.
	As an application we are able to effectively
	find finite q-product factorisations when they
	exist thus answering a question of Andrews.
	We provide other applications involving
	factorisations into theta-functions and eta-products.
\end{abstract}

\maketitle

\dedication{
Dedicated to George E. Andrews on the occasion of his 60th Birthday 
}

\tableofcontents

\section{Introduction}
\label{sec-intro}

	In the study of $q$-series one is quite often
interested in identifying generating functions as infinite
products. The classic example is the Rogers-Ramanujan
identity
$$
\sum_{n=0}^\infty \frac{q^{n^2}}{(q;q)_n}
=
\prod_{n=1}^\infty \frac{1}{(1-q^{5n-1})(1-q^{5n-4})}.
$$
Here we have used the notation in \eqn{qfact}.
It can be shown that the left-side of this identity
is the generating function for partitions whose parts
differ by at least two. The identity is equivalent
to saying the number of such partitions of $n$
is equinumerous with partitions of $n$ into parts congruent to
$\pm1\pmod{5}$.

	The main goals of the package are to provide
facility for handling the  following problems.
\begin{enumerate}
\item
Conversion of a given $q$-series into  an ``infinite'' product.
\item
Factorization of a given rational function into a finite $q$-product
if one exists.
\item
Find algebraic relations (if they exist) among the
$q$-series in a given list.
\end{enumerate}
	A $q$-product has the form
\begin{equation}
\prod_{j=1}^N (1-q^j)^{b_j},
\label{eq:qproddef}
\end{equation}
where the $b_j$ are integers.

	In  \cite[\S 10.7]{andrews:cbms}, George Andrews also
considered Problems 1 and 2, and asked for an easily accessible 
implementation. We provide implementations as well as considering
factorisations into theta-products and eta-products.
The package provides some basic functions for computing
$q$-series expansions of eta functions, theta functions,
Gaussian polynomials and $q$-products.
It also has a function for sifting out coefficients of
a $q$-series. It also has the basic infinite product
identities: the triple product identity, the quintuple
product identity and Winquist's identity.

\subsection{Installation instructions}
\label{sec-install}

The \t{qseries} package can be download via the WWW.
First use your favorite browser to access the url:\ssn
\t{http://www.math.ufl.edu/$\sim$frank/qmaple.html}
then follow the directions on that page. There are two versions:
one for UNIX and one for WINDOWS.

\section{Basic functions}
\label{sec-basfuncs}

	We describe the basic functions in the package which
are used to build $q$-series.
\subsection{Finite $q$-products}
\label{sec-qprods}
\subsubsection{Rising $q$-factorial}
\label{sec-qfac}
{\tt aqprod(a,q,n)} returns the product
\begin{equation}
(a;q)_n = (1-a)(1-aq)\cdots(1-aq^{n-1}).
\label{eq:qfact}
\end{equation}
We also use the notation
$$
(a;q)_\infty = \prod_{n=1}^\infty(1-aq^{n-1}).
$$

\subsubsection{Gaussian polynomials}
\label{sec-qbin}
When $0 \le m \le  n$,  \t{qbin(q,m,n)} returns the Gaussian 
polynomial (or $q$-binomial coefficient)
$$
\left [\begin {array}{c} n\\\noalign{\medskip}m\end {array}\right ]_q
={\frac { (q)_{{n}}}{ (q)_{{m}} (q)_{{n-m}}}},
$$
otherwise it returns $0$.

\subsection{Infinite products}
\label{sec-infprods}
\subsubsection{Dedekind eta products}
\label{sec-etaprods}
Suppose $\Re \tau>0$, and $q=\exp(2\pi i\tau)$. The Dedekind eta
function \cite[p. 121]{koblitz} is defined by
\begin{align*}
\eta(\tau) &= \exp(\pi i\tau/12)\prod_{n=1}^\infty
(1 -\exp(2\pi i n\tau))\\
           &= q^{1/24}\prod_{n=1}^\infty(1 -q^n).
\end{align*}
{\tt etaq(q,k,T)} returns the $q$-series expansion (up to $q^T$)
of the eta product
$$
\prod_{n=1}^\infty (1 - q^{kn}).
$$
This corresponds to the eta function $\eta(k\tau)$ except
for a power of $q$. Eta products occur frequently in the study
of $q$-series. For example, the generating function
for $p(n)$, the number of partitions of $n$, can be written
as
$$
\sum_{n=0}^\infty p(n) q^n = \frac{1}{\prod_{n=1}^\infty(1-q^n)}.
$$
See \cite[pp. 3--4]{andrews:book}. The generating function for the number of partitions
of $n$ that are $p$-cores \cite{gsk}, $a_p(n)$, can be written
$$
\sum_{n=0}^\infty a_p(n) q^n = \prod_{n=1}^\infty
\frac{(1-q^{pn})^p}{(1-q^n)}.
$$
Recently, Granville and Ono \cite{granvilleono} were able to prove a
long-standing conjecture in group representation theory
using elementary and function-theoretic properties of
the eta product above.

\subsubsection{Theta functions}
\label{sec-thetafuncs}
	Jacobi \cite[Vol I, pp. 497--538]{jacobi:works}
 defined four theta functions
$\theta_i(z,q)$, $i=1,2,3,4$. See also \cite[Ch. XXI]{ww:book}. Each theta function
can be written in terms of the others using a simple change
of variables. For this reason, it is common to define
$$
\theta(z,q) :=
\sum _{n=-\infty}^{\infty}{z}^{n}{q}^{{n}^{2}}.
$$
\t{theta(z,q,T)} returns the truncated theta-series                    
$$
\sum _{i=-T}^{T}{z}^{i}{q}^{{i}^{2}}.
$$
The case $z=1$ of Jacobi's theta functions occurs
quite frequently. We define
\begin{align*}
\theta_2(q) &:= \sum_{n=-\infty}^\infty q^{(n+1)^2/2} \\
\theta_3(q) &:= \sum_{n=-\infty}^\infty q^{n^2} \\
\theta_4(q) &:= \sum_{n=-\infty}^\infty (-1)^n q^{n^2} 
\end{align*}
\t{theta2(q,t)},
\t{theta3(q,t)},
\t{theta4(q,t)} (resp.)
returns the $q$-series expansion to
order $O\left(q^T\right)$ of          
$\theta_2(q)$,
$\theta_3(q)$, 
$\theta_4(q)$ respectively. 
Let $a$, and $b$ be positive integers and suppose $|q|<1$. 
Infinite products of the form
$$
(q^a; q^b)_\infty (q^{b-a}; q^b)_\infty
$$
occur quite frequently in the theory of partitions and $q$-series.
For example the right side of the Rogers-Ramanujan identity
is the reciprocal of the product with $(a,b)=(1,5)$.
The function {\tt jacprod(a,b,q,T)} returns the $q$-series expansion to 
order $O(q^T)$ of the  Jacobi-type infinite product given above.

\section{Product Conversion}
\label{sec-prodcon}

	In \cite[p. 233]{andrews:book}, \cite[\S 10.7]{andrews:cbms}
there is a very nice and useful algorithm
for converting a $q$-series into an infinite product. Any given  $q$-series
may be written formally as an infinite product
$$
1+ \sum_{n=1}^\infty b_n q^n = \prod_{n=1}^\infty (1-q^n)^{-a_n}.
$$
Here we assume the $b_n$ are integers. By taking the logarithmic
derivative of both sides we can obtain the recurrence
$$              
n b_n = \sum_{j=1}^n b_{n-j}\sum_{d\mid j} d a_d.
$$
Letting $a_n=1$ we obtain well-known special case
$$
n p(n) = 
 \sum_{j=1}^n p(n-j) \sigma(j).                
$$
We can also easily construct a recurrence for the $a_n$ from
the recurrence above.

  	The function \t{prodmake} is an implementation of Andrews'
algorithm. Other related functions are \t{etamake} and \t{jacprodmake}.

\subsection{\t{prodmake}}

	\t{prodmake(f,q,T)} converts the $q$-series $f$ into an infinite
product that agrees with $f$ to $\mbox{O}(q^T)$. Let's take a look at
the left side of the Roger's Ramanujan identity.
\sn
{\tt
$>$ \enskip with(qseries): \ssn
$>$ \enskip x:=1:  \ssn
$>$ \enskip for n from 1 to 8 do \ssn
\qquad\qquad x := x + q\carrot(n*n)/aqprod(q,q,n):\ssn
\hphantom{$>$} \enskip od:\ssn
$>$ \enskip x := series(x,q,50);\ssn
\begin{eqnarray*}
\lefteqn{x := 1 + q + q^{2} + q^{3} + 2\,q^{4} + 2\,q^{5} + 3\,q
^{6} + 3\,q^{7} + 4\,q^{8} + 5\,q^{9} + 6\,q^{10}} \\
 & &  + 7\,q^{11} + 9\,q^{12} + 10\,q^{13} + 12\,q^{14} + 14\,q^{
15} + 17\,q^{16} + 19\,q^{17} + 23\,q^{18} \\
 & &  + 26\,q^{19} + 31\,q^{20} + 35\,q^{21} + 41\,q^{22} + 46\,q^{23} + 54\,q^{24} + 61\,q^{25} + 70 \\
 & & q^{26} + 79\,q^{27} + 91\,q^{28} + 102\,q^{29} + 117\,q^{30} + 131\,q^{31} + 149\,q^{32} + 167 \\
 & & q^{33} + 189\,q^{34} + 211\,q^{35} + 239\,q^{36} + 266\,q^{37} + 299\,q^{38} + 333\,q^{39} +  \\
 & & 374\,q^{40} + 415\,q^{41} + 465\,q^{42} + 515\,q^{43} + 575\,q^{44} + 637\,q^{45} + 709\,q^{46} \\
 & &  + 783\,q^{47} + 871\,q^{48} + 961\,q^{49} + {\rm O}(q^{50})
\end{eqnarray*}
$>$ \enskip prodmake(x,q,40);
\begin{align*}
& 1/\left( 
(1-q)(1-{q}^{4})
(1-{q}^{6})(1-{q}^{9})(1-{q}^{11}
)(1-{q}^{14})(1-{q}^{16})
(1-{q}^{19})
\right.\\
& \left.
(1-{q}^{21})(1-
{q}^{24})(1-{q}^{26})(1-{q}^{29}
)(1-{q}^{31})(1-{q}^{34})
(1-{q}^{36})(1-{q}^{39})
\right)
\end{align*}
}
We have rediscovered the right side of the Rogers-Ramanujan
identity!

\EX{1}{Find (and prove) a product form for the $q$-series
$$
\sum_{n=0}^\infty \frac{q^{n^2}}{(q;q)_{2n}}.
$$
}
\subsection{\t{qfactor}}
The function \t{qfactor} is a version of \t{prodmake}. \break
\t{qfactor(f,T)} 
attempts
to write a rational function $f$ in $q$ as a $q$-product, ie. as
a product of terms of the form $(1-q^i)$. The second argument $T$
is optional. It specifies an an upper bound for the exponents of
$q$ that occur in the product. If $T$ is not specified it is given 
a default value of $4d+3$ where $d$ is the maximum of the degree
in $q$ of the numerator and denominator.  
The algorithm is quite simple. First the function is factored as usual,
and then it uses \t{prodmake} to do further factorisation into $q$-products.
Thus even if only part of the function can be written as a $q$-product
{\t qfactor} is able to find it.

As an example we consider some rational functions $T(r,h)$
introduced by Andrews \cite{andrews:cbms}(p.14) to explain
Rogers's \cite{rogers:1894} first proof of the Rogers-Ramanujan
identities. The $T(r,n)$ are defined recursively as follows:
\begin{align}
T(r,0) &= 1, \label{eq:rogpoly1}\\
T(r,1) &= 0, \label{eq:rogpoly2}\\
T(r,N) &= - \sum_{1\le 2j \le N} \left[
\begin{array}{c}
r + 2j\\
j
\end{array}\right]
T(r+2j,N-2j). \label{eq:rogpoly3}
\end{align}
\sn{\tt
$>$\ T:=proc(r,j) \sssn
$>$\ \ \ \ option\ remember; \sssn
$>$\ \ \ \ local\ x,k; \sssn
$>$\ \ \ \ x:=0; \sssn
$>$\ \ \ \ if\ j=0\ or\ j=1\ then \sssn
$>$\ \ \ \ \ \ RETURN((j-1)\carrot2): \sssn
$>$\ \ \ \ else \sssn
$>$\ \ \ \ \ \ \ for\ k\ from\ 1\ to\ floor(j/2)\ do \sssn
$>$\ \ \ \ \ \ \ \ \ \ \ x:=x-qbin(q,k,r+2*k)*T(r+2*k,j-2*k); \sssn
$>$\ \ \ \ \ \ \ \ od: \sssn
$>$\ \ \ \ \ \ \ \ RETURN(expand(x)); \sssn
$>$\ \ \ \ fi: \sssn
$>$\ end: \sssn
$>$ \enskip t8:=T(8,8);
\begin{eqnarray*}
\lefteqn{{\it t8} := 3\,q^{9} + 21\,q^{16} + 36\,q^{25} + 9\,q^{
36} + q^{6} + q^{7} + 2\,q^{8} + 5\,q^{10} + 6\,q^{11}} \\
 & & \mbox{} + 9\,q^{12} + 11\,q^{13} + 15\,q^{14} + 17\,q^{15}
 + 33\,q^{28} + 34\,q^{27} + 37\,q^{26} \\
 & & \mbox{} + 38\,q^{24} + 36\,q^{23} + 37\,q^{22} + 34\,q^{21}
 + 33\,q^{20} + 29\,q^{19} + 28\,q^{18} \\
 & & \mbox{} + 23\,q^{17} + 5\,q^{38} + 6\,q^{37} + 11\,q^{35} + 
15\,q^{34} + 17\,q^{33} + 21\,q^{32} \\
 & & \mbox{} + 23\,q^{31} + 28\,q^{30} + 29\,q^{29} + 3\,q^{39} + q^{42} + q^{41} + 2\,q^{40}\mbox{\hspace{62pt}}
\end{eqnarray*}
$>$ \enskip factor(t8);
\begin{eqnarray*}
\lefteqn{{q}^{6}\left ({q}^{4}+{q}^{3}+{q}^{2}+q+1\right )\left ({q}^{4}-{q}^{3
}+{q}^{2}-q+1\right )}\\
& & \mbox{ } \left ({q}^{10}+{q}^{9}+{q}^{8}+{q}^{7}+{q}^{6}+{
q}^{5}+{q}^{4}+{q}^{3}+{q}^{2}+q+1\right )\left ({q}^{4}+1\right )\\
& & \mbox{ } \left ({q}^{6}+{q}^{3}+1\right )\left ({q}^{8}+1\right )
\end{eqnarray*}
$>$ \enskip qfactor(t8,20);
\[
{\displaystyle \frac {(1 - q^{9})\,(1 - q^{10})\,(1 - q^{11})\,(1
 - q^{16})\,q^{6}}{(1 - q)\,(1 - q^{2})\,(1 - q^{3})\,(1 - q^{4})
}} 
\]
}
\sn
Observe how we used {\bf factor} to factor {\bf t8} into
cyclotomic polynomials. However, {\bf qfactor} was
able to factor {\bf t8} as a $q$-product.
We see that
$$
T(8,8) = \frac{(q^9;q)_3 (1-q^{16})q^6}{(q;q)_4}.
$$
\EX{2}{Use {\bf qfactor} to factorize $T(r,n)$ for different values
of $r$ and $n$. Then write $T(r,n)$ (defined above) as a $q$-product for general
$r$ and $n$.
}
For our next example we examine the sum 
$$
\sum_{k=-\infty}^\infty (-1)^k q^{k (3 k+1)/2}
\left[\begin{array}{c}
b+c\\
c+k
\end{array}\right]\,
\left[\begin{array}{c}
c+a\\
a+k
\end{array}\right]\,
\left[\begin{array}{c}
a+b\\
b+k
\end{array}\right]. 
$$
\sn{\tt
$>$\ dixson:=proc(a,b,c,q) \ssn
$>$\ \ \ \ local\ x,k,y; \ssn
$>$\ \ \ \ x:=0:\ \ y:=min(a,b,c): \ssn
$>$\ \ \ \ for\ k\ from\ -y\ to\ y\ do \ssn
$>$\ \ \ \ \ \ x:=x+(-1)\carrot{k}*q\carrot(k*(3*k+1)/2)*\ \ \ \ \  \ssn
$>$\ \ \ \ \ \ qbin(q,c+k,b+c)*qbin(q,a+k,c+a)*qbin(q,b+k,a+b);\ \  \ssn
$>$\ \ \ \ od: \ssn
$>$\ \ \ \ RETURN(x): \ssn
$>$\ end: \ssn
$>$ \enskip dx := expand(dixson(5,5,5,q)): \sn
$>$ \enskip qfactor(dx,20);
\begin{eqnarray*}
\lefteqn{(1 - q^{6})\,(1 - q^{7})\,(1 - q^{8})\,(1 - q^{9})\,(1
 - q^{10})\,(1 - q^{11})\,(1 - q^{12})} \\
 & & (1 - q^{13})\,(1 - q^{14})\,(1 - q^{15})/((1 - q)^{2}\,(1 - q^{2})^{2}\,(1 - q^{3})^{2}\,(1 - q^{4
})^{2} \\
 & & (1 - q^{5})^{2})
\end{eqnarray*}
}
We find that 
\begin{equation}
\sum_{k=-\infty}^\infty (-1)^k q^{k (3 k+1)/2}
\left[\begin{array}{c}
10\\
5+k
\end{array}\right]^3 = \frac{(q^6;q)_{10}}{(q;q)_5^2}.
\label{eq:dixeg}
\end{equation}
\EX{3}{Write the sum   
\[
\sum_{k=-\infty}^\infty (-1)^k q^{k (3 k+1)/2}
\left[\begin{array}{c}
2a\\
a+k
\end{array}\right]^3 
\]
as a $q$-product for general integral $a$.
}

\subsection{\t{etamake}}
Recall from  \sec{etaprods} that \t{etaq} is the function
to use for computing $q$-expansions of eta products. 
If one wants to apply the theory of modular forms to $q$-series
it is quite useful to determine whether a given $q$-series
is a product of eta functions. The function in the package
for doing this conversion is \t{etamake}.
\t{etamake(f,q,T)} will write the given $q$-series $f$ as a product
of eta functions which agrees with $f$ up to $q^T$.
As an example, let's see how we can write the theta functions
as eta products.
\sn{
\tt
$>$ \enskip theta2(q,100)/q\carrot(1/4);
\begin{eqnarray*}
\lefteqn{2\,{q}^{132} + 2\,{q}^{110} + 2\,{q}^{90} + 2\,{q}^{72}
 + 2\,{q}^{56} + 2\,{q}^{42} + 2\,{q}^{30} + 2\,{q}^{20} + 2\,{q}
^{12} + 2\,{q}^{6}} \\
 & & \mbox{} + 2\,{q}^{2} + 2 + {q}^{156}\mbox{\hspace{237pt}}
\end{eqnarray*}
$>$ \enskip etamake(",q,100);
\[
2\,{\displaystyle \frac {{ \eta}(\,4\,{ \tau}\,)^{2}}{{q}^{1/4}\,
{ \eta}(\,2\,{ \tau}\,)}}
\]
$>$ \enskip theta3(q,100);
\begin{eqnarray*}
\lefteqn{2\,{q}^{121} + 2\,{q}^{100} + 2\,{q}^{81} + 2\,{q}^{64}
 + 2\,{q}^{49} + 2\,{q}^{36} + 2\,{q}^{25} + 2\,{q}^{16} + 2\,{q}
^{9} + 2\,{q}^{4}} \\
 & & \mbox{} + 2\,{q} + 1\mbox{\hspace{268pt}}
\end{eqnarray*}
$>$ \enskip etamake(",q,100);
\[
{\displaystyle \frac {{ \eta}(\,2\,{ \tau}\,)^{5}}{{ \eta}(\,4\,{
 \tau}\,)^{2}\,{ \eta}(\,{ \tau}\,)^{2}}}
\]
$>$ \enskip theta4(q,100);
\begin{eqnarray*}
\lefteqn{ - 2\,{q}^{121} + 2\,{q}^{100} - 2\,{q}^{81} + 2\,{q}^{
64} - 2\,{q}^{49} + 2\,{q}^{36} - 2\,{q}^{25} + 2\,{q}^{16} - 2\,
{q}^{9} + 2\,{q}^{4}} \\
 & & \mbox{} - 2\,{q} + 1\mbox{\hspace{264pt}}
\end{eqnarray*}
$>$ \enskip etamake(",q,100);
\[
{\displaystyle \frac {{ \eta}(\,{ \tau}\,)^{2}}{{ \eta}(\,2\,{ 
\tau}\,)}}
\]
}
We are led to the well-known identities:
\begin{align*}
\theta_2(q) &= 2 \frac{\eta(4\tau)^2}{\eta(2\tau)},\\
\theta_3(q) &=  \frac{\eta(2\tau)^5}{\eta(4\tau)^2 \eta(\tau)^2},\\
\theta_4(q) &=  \frac{\eta(\tau)^2}{\eta(2\tau)}.
\end{align*}
The idea of the algorithm is quite simple. Given a $q$-series $f$
(say with leading coefficient 1) one just keeps recursively
multiplying by powers of the right eta function until  the desired  
terms agree. For example, suppose we are given a $q$-series
$$
f = 1 + c q^k + \cdots.
$$
Then the next step is to multiply by \t{etaq(q,k,T)\carrot(-c)}.
\EX{4}{
Define the $q$-series
\begin{align*}
a(q) &:= \sum_{n=-\infty}^\infty \sum_{m=-\infty}^\infty
         q^{n^2+nm+m^2} \\
b(q) &:= \sum_{n=-\infty}^\infty \sum_{m=-\infty}^\infty
         \omega^{n-m} q^{n^2+nm+m^2} \\
c(q) &:= \sum_{n=-\infty}^\infty \sum_{m=-\infty}^\infty
         q^{(n+1/3)^2+(n+1/3)(m+1/3)+(m+1/3)^2} 
\end{align*}
where $\omega=\exp(2\pi i/3)$.
Two of the three functions above can be written as eta products.
Can you find them? \\
{\it Hint}\,: It would be wise to define
\sn{
\tt
$>$ \enskip omega := RootOf(z\carrot2 + z + 1 = 0);
}}
\subsection{\t{jacprodmake}}
In \sec{thetafuncs} we observed that the right side of
the Rogers-Ramanujan identity could be a written in terms
of a Jacobi product. The function \t{jacprodmake} 
converts
a $q$-series into a Jacobi-type product if one exists.
Given a $q$-series $f$, 
\t{jacprodmake(f,q,T)} attempts to convert  $f$ into a product of 
theta functions that agrees with
$f$ to order $\mbox{O}(q^T)$.
Each theta-function has the form $\mbox{JAC}(a,b,\infty)$, where 
$a$, $b$ are integers and $0 \le a < b$.
If $0 < a$, then $\mbox{JAC}(a,b,\infty)$ corresponds to the 
theta-product
$$
(q^a;q^b)_\infty(q^{b-a};q^b)_\infty(q^b;q^b)_\infty.
$$
We call this a theta product because it is $\theta(-q^{(b-2a)/2},q^{b/2})$.
The \t{jacprodmake} function is really a variant
of \t{prodmake}. It involves using \t{prodmake} to compute
the sequence of exponents and then searching for periodicity.

   If $a = 0$, then $\mbox{JAC}(0,b,\infty)$ corresponds to the 
eta-product
$$
                             (q^b;q^b)_\infty.
$$
Let's re-examine the Rogers-Ramanujan identity.
\sn{\tt
$>$ \enskip with(qseries): \ssn
$>$ \enskip x:=1: \ssn
$>$ \enskip for n from 1 to 8 do \ssn
$>$ \qquad    x:=x+q\carrot(n*n)/aqprod(q,q,n): \ssn
$>$ \enskip od: \ssn
$>$ \enskip x:=series(x,q,50): \ssn
$>$ \enskip y:=jacprodmake(x,q,40);
$$
y := {\frac {{\it JAC}(0,5,\infty )}{{\it JAC}(1,5,\infty )}}
$$
$>$ \enskip z:=jac2prod(y);
$$
z := {\frac {1}{ (q,{q}^{5})_{{\infty }} ({q}^{4},{q}^{5})_{{\infty }}}}
$$
}
Note that we were able to observe that the left side
of the Rogers-Ramanujan identity (at least up through $q^{40}$)
can be written as a quotient of theta functions. We used
the function \t{jac2prod}, to simplify the result and get it into
a more recognizable form. The function 
\t{jac2prod(jacexpr)} converts a product of theta functions into
$q$-product form; ie. as a product of 
functions of the form  $(q^a ; q^b)_\infty$.
Here \t{jacexpr}  is a product (or quotient)
of terms $\mbox{JAC}(i,j, \infty)$, where $i$, $j$ are integers
and $0 \le i < j$.

A related function is \t{jac2series}. This converts a Jacobi-type
product into a form better for computing its $q$-series.
It simply replaces each Jacobi-type product with its
corresponding theta-series.
\sn{\tt
$>$ \enskip with(qseries):\ssn
$>$ \enskip x:=0:\ssn
$>$ for n from 0 to 10 do \ssn
\qquad x := x + q\carrot(n*(n+1)/2)*aqprod(-q,q,n)/aqprod(q,q,2*n+1):\ssn
\quad od: \ssn
$>$ \enskip x:=series(x,q,50):\ssn
$>$ \enskip jacprodmake(x,q,50);
\begin{eqnarray*}
\lefteqn{{\rm JAC}(0, \,14, \,\infty )^{6} \left/ {\vrule 
height0.44em width0em depth0.44em} \right. \!  \! ({\rm JAC}(1, 
\,14, \,\infty )^{2}\,{\rm JAC}(3, \,14, \,\infty )\,{\rm JAC}(4
, \,14, \,\infty )} \\
 & & {\rm JAC}(5, \,14, \,\infty )\,{\rm JAC}(6, \,14, \,\infty )
\,\sqrt{{\displaystyle \frac {{\rm JAC}(7, \,14, \,\infty )}{
{\rm JAC}(0, \,14, \,\infty )}} })\mbox{\hspace{74pt}}
\end{eqnarray*}
$>$ \enskip jac2series(",500);
\begin{eqnarray*}
\lefteqn{(q^{364} - q^{210} + q^{98} - q^{28} + 1 - q^{14} + q^{
70} - q^{168} + q^{308} - q^{490})^{6} \left/ {\vrule height0.44em width0em depth0.44em}
 \right. \!  \! ((} \\
 & &  - q^{621} + q^{496} - q^{385} + q^{288} - q^{205} + q^{136}
 - q^{81} + q^{40} - q^{13} + 1 - q \\
 & & \mbox{} + q^{16} - q^{45} + q^{88} - q^{145} + q^{216} - q^{
301} + q^{400} - q^{513})^{2}( - q^{603} + q^{480} \\
 & & \mbox{} - q^{371} + q^{276} - q^{195} + q^{128} - q^{75} + q
^{36} - q^{11} + 1 - q^{3} + q^{20} - q^{51} \\
 & & \mbox{} + q^{96} - q^{155} + q^{228} - q^{315} + q^{416} - q
^{531})( - q^{594} + q^{472} - q^{364} + q^{270} \\
 & & \mbox{} - q^{190} + q^{124} - q^{72} + q^{34} - q^{10} + 1 - q^{4} + q^{22} - q^{54} + q^{100} - q^{160} \\
 & & \mbox{} + q^{234} - q^{322} + q^{424} - q^{540})( - q^{585} + q^{464} - q^{357} + q^{264} - q^{185} + q^{120} \\
 & & \mbox{} - q^{69} + q^{32} - q^{9} + 1 - q^{5} + q^{24} - q^{57} + q^{104} - q^{165} + q^{240} - q^{329} \\
 & & \mbox{} + q^{432} - q^{549})( - q^{576} + q^{456} - q^{350} + q^{258} - q^{180} + q^{116} - q^{66} + q^{30} \\
 & & \mbox{} - q^{8} + 1 - q^{6} + q^{26} - q^{60} + q^{108} - q^{170} + q^{246} - q^{336} + q^{440} - q^{558})(( \\
 & &  - 2\,q^{567} + 2\,q^{448} - 2\,q^{343} + 2\,q^{252} - 2\,q^{175} + 2\,q^{112} - 2\,q^{63} + 2\,q^{28} \\
 & & \mbox{} - 2\,q^{7} + 1)/( \\
 & & q^{364} - q^{210} + q^{98} - q^{28} + 1 - q^{14} + q^{70} - 
q^{168} + q^{308} - q^{490}))^{1/2})
\end{eqnarray*}
}
It seems that the $q$-series
$$
\sum_{n\ge0} \frac{(-q;q)_n q^{n(n+1)/2}}{(q;q)_{2n+1}}
$$
can be written as Jacobi-type product. Assuming that this is the case
we used \t{jac2series} to write this $q$-series in terms of theta-series
at least up to $q^{500}$. This should provide an efficient method
for computing the $q$-series expansion and also for computing
the function at particular values of $q$.

\EX{5}{Use \t{jacprodmake} and \t{jac2series} to compute the
$q$-series expansion of
$$
\sum_{n\ge0} \frac{(-q;q)_n q^{n(3n+1)/2}}{(q;q)_{2n+1}}
$$
up to $q^{1000}$, assuming it is Jacobi-type product. Can you identify
the infinite product?
}

\section{The Search for Relations}
\label{sec-srels}

The functions for finding relations between $q$-series are
\t{findhom},\break
\t{findhomcombo}, \t{findnonhom}, 
\t{findnonhomcombo}, and \t{findpoly}.

\subsection{\t{findhom}}
\label{sec-findhom}

If the $q$-series one is concerned with are modular forms of
a particular weight, then theoretically these functions will 
satisfy homogeneous polynomial relations. See \cite[p. 263]{garvanrad},
\cite{bbg2} for more details and examples.
The function \t{findhom(L,q,n,topshift)}
returns a set of potential homogeneous relations
of order $n$ among the $q$-series in the list $L$.
The value of {\tt topshift} is usually taken to be zero. However if
it appears that spurious relations are being generated then a higher
value of \t{topshift} should be taken.

The idea is to convert this into a linear algebra problem.
This program generates a list of monomials of degree $n$ of the functions
in the given list of $q$-series $L$. 
The $q$-expansion (up to a certain point) of each monomial is found
and converted into a row vector of a matrix.
The set of relations is then found by computing
the kernel of the transpose of this matrix. As an example, we now consider 
relations between the theta functions $\theta_3(q)$, $\theta_4(q)$,
$\theta_3(q^2)$, and $\theta_4(q^2)$.

\sn
{\tt
$>$ \enskip  with(qseries):\ssn
$>$ \enskip  findhom([theta3(q,100),theta4(q,100),theta3(q\carrot2,100), 
\sssn
theta4(q\carrot2,100)],q,1,0);
\begin{center}
\# of terms , 25
\end{center}
\begin{center}
-----RELATIONS-----of order---, 1
\end{center}
$$
\left \{\left \{\right \}\right \}
$$
$>$ \enskip  findhom([theta3(q,100),theta4(q,100),theta3(q\carrot2,100), 
\sssn
theta4(q\carrot2,100)],q,2,0);
\begin{center}
\# of terms , 31
\end{center}
\begin{center}
-----RELATIONS-----of order---, 2
\end{center}
$$               
\left \{{X_{{1}}}^{2}  + {X_{{2}}}^{2}  - 2\,{X_{{3}}}^{2},  - X_{{1}}X_{{2}}  + {X_{{4}}}^{2}\right \} 
$$                
}\noindent
From the session above we see that there is no linear relation between
the functions
$\theta_3(q)$, $\theta_4(q)$, $\theta_3(q^2)$ and $\theta_4(q^2)$.
However, it appears that there are two quadratic relations:
$$
\theta_3(q^2) =\left(\frac{\theta_3^2(q)  + \theta_4^2(q)}{2}\right)^{1/2} 
$$
and
$$
\theta_4(q^2) = \left(\theta_3^2(q) \theta_4^2(q)\right)^{1/2}.
$$
This is Gauss' parametrization of the arithmetic-geometric
mean iteration. See \cite[Ch. 2]{borwein:book} for details.

\EX{6}{
Define $a(q)$, $b(q)$ and $c(q)$ as in {\bf Exercise 2}.
Find homogeneous relations between the functions
$a(q)$, $b(q)$, $c(q)$, 
$a(q^3)$, $b(q^3)$, $c(q^3)$. In particular, try to get
$a(q^3)$ and  $b(q^3)$ in terms of $a(q)$, $b(q)$.
See \cite{bbg} for more details.  These results lead to
a cubic analog of the AGM due to Jon and Peter Borwein 
\cite{bbcubic}, \cite{bbcubic2}.
}

\subsection{\t{findhomcombo}}
\label{sec-findhomcombo}

The function \t{findhomcombo} is a variant of \t{findhom}.
Suppose $f$ is a $q$-series and $L$ is a list of $q$-series.
\t{findhomcombo(f,L,q,}\break
\t{n,topshift,etaoption)} tries to express $f$
as a homogeneous polynomial in the members of $L$.
If \t{etaoption=yes} then each monomial in the 
combination is ``converted''
into an eta-product using \t{etamake}.

We illustrate this function with certain Eisenstein series.
For $p$ an odd prime define 
$$
\chi(m)=\left(\frac{m}{p}\right)\qquad\mbox{(the Legendre symbol)}.
$$
Suppose $k$ is an  integer, $k\ge2$, and $(p-1)/2\equiv k\pmod{2}$.
Define the Eisenstein series
$$
U_{p,k}(q) := \sum_{m=1}^\infty\sum_{n=1}^\infty \chi(m) n^{k-1}q^{mn}.
$$
Then $U_{p,k}$ is a modular form of weight $k$ and character $\chi$ for
the congruence subgroup $\Gamma_0(p)$. See \cite{kolberg},
\cite{garvantcore} 
for more details. The classical result is the following identity
found by Ramanujan \cite[Eq. (1.52), p. 354]{ramlost}:
$$
U_{5,2} = \frac{\eta(5\tau)^5}{\eta(\tau)}.
$$
Kolberg \cite{kolberg} has found many relations
between such Eisenstein series and certain eta products.
The eta function $\eta(\tau)$ is a modular form of weight
$\tfrac{1}{2}$ \cite[p. 121]{koblitz}. Hence the modular forms
$$
B_1:=\frac{\eta(5\tau)^5}{\eta(\tau)},\qquad
B_2:=\frac{\eta(\tau)^5}{\eta(5\tau)}       
$$
are modular forms of weight $\tfrac{(5-1)}{2}=2$.
In fact, it can be shown that they are modular forms on $\Gamma_0(5)$
with character $\left(\frac{\cdot}{5}\right)$.
We might therefore expect that $U_{5,6}$ can be written
as a homogeneous cubic polynomial in $B_1$, $B_2$.
We write a short \maple program to compute the Eisenstein  
series $U_{p,k}$.
\sn{\tt
$>$ \enskip with(numtheory): \ssn
$>$ \enskip UE:=proc(q,k,p,trunk) \ssn
$>$ \enskip \enskip  local x,m,n: \ssn
$>$ \enskip \enskip  x:=0: \ssn
$>$ \enskip \enskip  for m from 1 to trunk do \ssn
$>$ \enskip \enskip  \enskip   for n from 1 to trunk/m do  \ssn
$>$ \enskip \enskip  \enskip \enskip x:=x + L(m,p)*n\carrot(k-1)*q\carrot(m*n):  \ssn
$>$ \enskip \enskip  \enskip   od: \ssn
$>$ \enskip \enskip  od:  \ssn
$>$ \enskip end: \ssn
}
The function \t{UE(q,k,p,trunk)} returns the $q$-expansion
of $U_{p,k}$ up through $q^{\mbox{trunk}}$.
We note that \t{L(m,p)} returns the Legendre symbol
$\displaystyle \left(\frac{m}{p}\right)$.
We are now ready to study $U_{5,6}$.
\sn
{\tt
$>$ \enskip with(qseries): \ssn
$>$ \enskip f := UE(q,6,5,50): \ssn
$>$ \enskip B1 := etaq(q,1,50){\carrot}5/etaq(q,5,50): \ssn
$>$ \enskip B2 := q*etaq(q,5,50){\carrot}5/etaq(q,1,50): \ssn
$>$ \enskip findhomcombo(f,[B1,B2],q,3,0,yes);
\begin{center}
\# of terms , 25
\end{center}
\begin{center}
-----possible linear combinations of degree------, 3
\end{center}
$$
\left \{\eta(5\,\tau)^{3}\eta(\tau)^{9}+40\,\eta(5\,\tau)^{9}
\eta(\tau)^{3}+335\,{\frac {\eta(5\,\tau)^{15}}{\eta(\tau)^{3}}}\right \}
$$
$$
\left \{{X_{{1}}}^{2}X_{{2}} 
+ 40\,X_{{1}}{X_{{2}}}^{2} 
+ 335\,{X_{{2}}}^{
3}\right \}
$$
}
It would appear that                                      

$$
U_{5,6} = \eta(5\,\tau)^{3}\eta(\tau)^{9}+40\,\eta(5\,\tau)^{9}
\eta(\tau)^{3}+335\,{\frac {\eta(5\,\tau)^{15}}{\eta(\tau)^{3}}}.
$$
The proof is a straightforward exercise using the theory of
modular forms.
\EX{7}{
Define the following eta products:
$$
C_1:=\frac{\eta(7\tau)^7}{\eta(\tau)},\qquad
C_2:=\eta(\tau)^3\,\eta(7\tau)^3,\qquad
C_3:=\frac{\eta(\tau)^7}{\eta(7\tau)}.      
$$
What is the weight of these modular forms? \sn
Write $U_{7,3}$ in terms of $C_1$, $C_2$, $C_3$.\sn
The identity that you should find was originally due to   
Ramanujan.
Also see Fine \cite[p. 159]{fine} and \cite[Eq. (5.4)]{gsk}.\sn
If you are ambitious find $U_{7,9}$ in terms
of $C_1$, $C_2$, $C_3$.\sn
}

\subsection{\t{findnonhom}}
\label{sec-findnonhom}

In section \ref{sec-findhom} we introduced the function
\t{findhom} to find homogeneous relations between $q$-series.
The nonhomogeneous analog is
\t{findnonhom}.

The syntax of \t{findnonhom} is the same as \t{findhom}.
Typically (but not necessarily) \t{findhom} is used to find relations
between modular forms of a certain weight. To find relations
between modular functions we would use \t{findnonhom}.
We consider an example involving theta functions.
\sn
{\tt
$>$ \enskip with(qseries): \ssn
$>$ \enskip F := q -> theta3(q,500)/theta3(q\carrot5,100): \ssn
$>$ \enskip U := 2*q*theta(q\carrot10,q\carrot25,5)/theta3(q\carrot25,20);
$$
U:={\frac {2\,q\left ({q}^{575}+{q}^{360}+{q}^{195}+{q}^{80}+{q}^{15}+1
+{q}^{35}+{q}^{120}+{q}^{255}+{q}^{440}+{q}^{675}\right )}{2\,{q}^{625
}+2\,{q}^{400}+2\,{q}^{225}+2\,{q}^{100}+2\,{q}^{25}+1}}
$$
$>$ \enskip EQNS := findnonhom([F(q),F(q\carrot5),U],q,3,20);
\begin{center}
\# of terms , 61
\end{center}
\begin{center}
-----RELATIONS-----of order---, 3
\end{center}
$$             
 \mbox {{\tt EQNS}}:=\left \{  - 1  - X_{{1}}X_{{2}}X_{{3}}  + {X_{{2}}}^{2}  + {X_{{3}}}^{2}  + X_{{3}}\right \} 
$$            
$>$ \enskip ANS:=EQNS[1];
$$
\mbox{ANS}:=-1-X_{{1}}X_{{2}}X_{{3}}+{X_{{2}}}^{2}+{X_{{3}}}^{2}+X_{{3}}
$$
$>$ \enskip CHECK := subs({X[1]=F(q),X[2]=F(q\carrot5),X[3]=U},ANS): \ssn
$>$ \enskip series(CHECK,q,500);
$$
O\left ({q}^{500}\right )
$$
}
\sn

We define
$$
F(q) := \frac{\theta_3(q)}{\theta_3(q^5)},     
$$
and
$$
U(q)
:= 2\displaystyle\frac{\displaystyle\sum_{n=-\infty}^\infty q^{25n^2+10n+1}}
{\theta_3(q^{25})}.
$$
We note that $U(q)$ and $F(q)$ are modular functions since they
are ratios of theta series.
From the session above we see that it appears that 
$$
1+F(q)F(q^5)U(q)={F(q^5)}^{2}+{U(q)}^{2}+U(q).
$$
Observe how we were able to 
verify  this equation to high order.
When \t{findnonhom} returns a set of relations the variable $X$
has been declared {\it global}. This is so we can manipulate
the relations. It this way we were able to assign \t{ANS}
to the relation found and then use \t{subs} and \t{series}
to check it to order $\mbox{O}\left(q^{500}\right)$.

\subsection{\t{findnonhomcombo}}
\label{sec-findnonhomcombo}

The syntax of \t{findnonhomcombo} is the same as 
\t{findhomcombo}.
We consider an example involving eta functions.
First we define 
$$
\xi:={\frac {\eta(49\,\tau)}{\eta(\tau)}},
$$
and
$$
T:=\left (\frac{\eta(7\,\tau)}{\eta(\tau)}\right)^{4},
$$
Using the theory of modular functions it can be shown
that one must be able to write $T^2$ in terms of $T$
and powers of $\xi$. 
We now use \t{findnonhomcombo}
to get $T^2$ in terms of $\xi$ and $T$.
\sn
{\tt
$>$ \enskip  xi:=q{\carrot}2*etaq(q,49,100)/etaq(q,1,100):\ssn
$>$ \enskip  T:=q*(etaq(q,7,100)/etaq(q,1,100)){\carrot}4:\ssn
$>$ \enskip  findnonhomcombo(T{\carrot}2,[T,xi],q,7,-15);
\begin{center}
\# of terms , 42

matrix is , 37, x, 42

-----possible linear combinations of degree------, 7
\end{center}
\begin{align*}   
&\left \{147\,{X_{{2}}}^{5}  
+ 343\,{X_{{2}}}^{7}  + 343\,{X_{{2}}}^{6}  + X_{{2}}  
+ 49\,{X_{{2}}}^{4}  
+ 49\,X_{{1}}{X_{{2}}}^{3}  + 7\,X_{{1}}X_{{2}} \right. \\
&\qquad \left. + 21\,{X_{{2}}}^{3}  
+ 7\,{X_{{2}}}^{2} 
+ 35\,X_{{1}}{X_{{2}}}^{2}\right \}
\end{align*}          
}
\sn
Then it seems that
$$
{T}^{2}=\left (35\,{\xi}^{2}+49\,{\xi}^{3}+7\,\xi\right )T+
343\,{\xi}^{7}+343\,{\xi}^{6}+147\,{\xi}^{5}+49\,{\xi}^{4}+
21\,{\xi}^{3}+7\,{\xi}^{2}+\xi.
$$
This is the modular equation used by Watson\cite{watson38} to prove
Ramanujan's partition congruences for powers of $7$.
Also see \cite{atkin} and \cite{knopp}, and see
\cite{garvan84} for an elementary treatment.

\EX{8}{
Define
$$
\xi:={\frac {\eta(25\,\tau)}{\eta(\tau)}},
$$
and
$$
T:=\left (\frac{\eta(5\,\tau)}{\eta(\tau)}\right)^{6}.
$$
Use \t{findnonhomcombo} to express $T$ as a polynomial
in $\xi$ of degree $5$. The modular equation you find
was used by Watson to prove Ramanujan's partition congruences
for powers of $5$. See \cite{hirschhorn81} for an elementary
treatment.
}

\EX{9}{Define $a(q)$ and $c(q)$ as in {\bf Exercise 2}.
Define 
$$
x(q):={\frac {{c(q)}^{3}}{{a(q)}^{3}}},
$$
and the classical Eisenstein series (usually called $E_3$; see 
\cite[p. 93]{serre})
$$
N(q):=1-504\,\sum _{n=1}^{\infty }{\frac {{n}^{5}{q}^{n}}{1-{q}^{n}}}.
$$
Use \t{findnonhomcombo} to express $N(q)$ in terms of $a(q)$ and $x(q)$.
\ssn
{\it HINT}\,: $N(q)$ is a modular form of weight $6$ and $a(q)$ and $c(q)$
are modular forms of weight $1$. 
See \cite{berndtbhargavagarvan} for this result and many others.
}

\subsection{\t{findpoly}}
\label{sec-findpoly}

The function \t{findpoly}
is used to find a polynomial relation between
two given $q$-series with degrees specified.
 
\t{findpoly(x,y,q,deg1,deg2,check)} returns a possible polynomial in 
$X$, $Y$ (with corresponding degrees deg1, deg2) which is satisfied by 
the $q$-series $x$ and $y$.
                                                              
If check is assigned then the relation is checked to 
$\mbox{O}\left(q^{\mbox{check}}\right)$.
 
We illustrate this function with an example involving theta functions
and the function $a(q)$ and $c(q)$ encountered in {\bf Exercises 2} and
{\bf 7}. It can be shown that
$$
a(q) = 
\theta_{{3}}(q)\theta_{{3}}({q}^{3})+\theta_{{2}}(q)\theta_{{2}}({q}
^{3}).
$$
See \cite{bbg} for details. This equation provides
a better way of computing the $q$-series expansion of $a(q)$ than
the definition. 
In {\bf Exercise 2} you would have found that
$$
c(q) = 3\,\
\frac{\eta^3(3\,\tau)}{\eta(\tau)}.
$$
See \cite{bbg} for a proof.
Define 
$$
y:={\frac {{c}^{3}}{{a}^{3}}},
$$
and  
$$
x:=
\left(\frac{\theta_{{2}}(q)}{\theta_{{2}}({q}^{3})}\right)^{2}
+
\left(\frac{\theta_{{3}}(q)}{\theta_{{3}}({q}^{3})}\right)^{2}.
$$
We use \t{findpoly} to find a polynomial
relation between $x$ and $y$.
\sn
{\tt
$>$ \enskip with(qseries):\ssn
$>$ \enskip x1 := radsimp(theta2(q,100){\carrot}2/theta2(q{\carrot}3,40){\carrot}2): \ssn
$>$ \enskip x2 := theta3(q,100){\carrot}2/theta3(q{\carrot}3,40){\carrot}2: \ssn
$>$ \enskip x := x1+x2: \ssn
$>$ \enskip c := q*etaq(q,3,100){\carrot}9/etaq(q,1,100){\carrot}3: \ssn
$>$ \enskip a := radsimp(theta3(q,100)*theta3(q{\carrot}3,40)+theta2(q,100)\ssn
\quad *theta2(q{\carrot}3,40)): \ssn
$>$ \enskip c := 3*q{\carrot}(1/3)*etaq(q,3,100){\carrot}3/etaq(q,1,100): \ssn
$>$ \enskip y := radsimp(c{\carrot}3/a{\carrot}3): \ssn
$>$ \enskip P1:=findpoly(x,y,q,3,1,60);\ssn
WARNING: X,Y are global.\ssn
\begin{center}
                                  dims , 8, 18

                               The polynomial is
\end{center}
$$
\left (X+6\right )^{3}Y-27\,\left (X+2\right )^{2}
$$
\begin{center}
                             Checking to order, 60
\end{center}
$$
\mbox{O}\left(q^{59}\right)
$$
$$
\mbox{P1}:=\left (X+6\right )^{3}Y-27\,\left (X+2\right )^{2}
$$
}
It seems that $x$ and $y$ satisfy the equation
$$
p(x,y)=(x+6)^3\,y - 27(x+2)^2=0.
$$
Therefore it would seem that 
$$
{\frac {{c}^{3}}{{a}^{3}}}=27\,{\frac {\left (x+2\right )^{2}}{\left (
x+6\right )^{3}}}.
$$
See \cite[pp. 4237--4240]{berndtbhargavagarvan} for more details.

\EX{10}{Define 
$$
m:=
\left(\frac{\theta_{{3}}(q)}{\theta_{{3}}({q}^{3})}\right)^{2}.
$$
Use \t{polyfind} to find $y=\frac{c^3}{a^3}$ as a rational function in $m$.
}

\section{Sifting coefficients}
\label{sec-sift}

Suppose we are given a $q$-series 
$$
A(q) = \sum_{n=0}^\infty a_n  q^n.
$$
Occasionally it will turn out the generating function
$$
\sum_{n=0}^\infty a_{mn+r} q^n
$$
will have a very nice form. A famous example for $p(n)$
is due to Ramanujan:
$$
\sum_{n=0}^\infty p(5n+4) q^n = 5 \prod_{n=1}^\infty \frac{(1-q^{5n})^6}
{(1-q^n)^5}.
$$
See \cite[Cor. 10.6]{andrews:book}. 
In fact, G.H. Hardy and Major MacMahon \cite[p. xxxv]{ramanujan:papers} both
agreed that this is Ramanujan's most beautiful identity.

Suppose $s$ is the $q$-series
$$
 \sum a_{i} q^i + \mbox{O}(q^T)
$$
then {\tt sift(s,q,n,k,T)} returns the $q$-series
$$
 \sum a_{ni+k} q^i
 + \mbox{O}(q^{T/n}).
$$

We illustrate this function with another example from
the theory of partitions. Let $pd(n)$ denote
the number of partitions of $n$ into distinct parts.
Then it is well known that
$$
\sum_{n=0}^\infty pd(n) q^n = \prod_{n=1}^\infty (1+q^n)
=\prod_{n=1}^\infty \frac{(1-q^{2n})}{(1-q^n)}.
$$
We now examine the generating function of $pd(5n+1)$ in \maplep
\sn{\tt
$>$ \enskip PD:=series(etaq(q,2,200)/etaq(q,1,200),q,200): \ssn
$>$ \enskip sift(PD,q,5,1,199);
\begin{eqnarray*}
\lefteqn{1 + 4\,{q} + 5010688\,{q}^{26} + 53250\,{q}^{15} + 668\,
{q}^{7} + 12\,{q}^{2} + 165\,{q}^{5}} \\
 & & \mbox{} + 12076\,{q}^{12} + 1087744\,{q}^{22} + 109420549\,{
q}^{35} + 76\,{q}^{4} + 32\,{q}^{3} \\
 & & \mbox{} + 340\,{q}^{6} + 1260\,{q}^{8} + 2304\,{q}^{9} + 
4097\,{q}^{10} + 7108\,{q}^{11} + 20132\,{q}^{13} \\
 & & \mbox{} + 32992\,{q}^{14} + 84756\,{q}^{16} + 133184\,{q}^{
17} + 206848\,{q}^{18} + 317788\,{q}^{19} \\
 & & \mbox{} + 728260\,{q}^{21} + 20792120\,{q}^{30} + 2368800\,{
q}^{24} + 483330\,{q}^{20} \\
 & & \mbox{} + 1611388\,{q}^{23} + 3457027\,{q}^{25} + 7215644\,{
q}^{27} + 10327156\,{q}^{28} \\
 & & \mbox{} + 14694244\,{q}^{29} + 29264960\,{q}^{31} + 40982540
\,{q}^{32} + 57114844\,{q}^{33} \\
 & & \mbox{} + 79229676\,{q}^{34} + 150473568\,{q}^{36} + 
206084096\,{q}^{37} \\
 & & \mbox{} + 281138048\,{q}^{38} + 382075868\,{q}^{39}
\end{eqnarray*}
$>$ \enskip PD1:=": \ssn
$>$ \enskip etamake(PD1,q,38);
\[
{\displaystyle \frac {{ \eta}(\,5\,{ \tau}\,)^{3}\,{ \eta}(\,2\,{
 \tau}\,)^{2}}{{q}^{5/24}\,{ \eta}(\,10\,{ \tau}\,)\,{ \eta}(\,{ 
\tau}\,)^{4}}}
\]
}

So it would seem that
$$
\sum_{n=0}^\infty pd(5n+1) q^n = \prod_{n=1}^\infty
\frac{ (1-q^{5n})^3 (1-q^{2n})^2}
{(1-q^{10n}) (1-q^n)^4}.
$$
This result was found originally by R{\o}dseth \cite{rodseth69}.

\EX{11}{R{\o}dseth also found the generating functions for
$pd(5n+r)$ for $r=0$, $1$, $2$, $3$ and $4$. For each $r$
use \t{sift} and \t{jacprodmake} to identify these
generating functions as infinite products.
}

\section{Product Identities}
\label{sec-prodids}

At present, the package contains the Triple Product identity,
the Quintuple Product identity and Winquist's identity.
These are the most commonly used of the Macdonald identities
\cite{macdonald72}, \cite{stanton86}, \cite{stanton88}. 
The Macdonald identities 
are the analogs of the Weyl denominator for affine roots systems.
Hopefully, a later version of this package will include these
more general identities.

\subsection{The Triple Product Identity}
\label{sec-tripleprod}

The triple product identity is
\begin{equation}
\sum_{n=-\infty}^\infty
(-1)^n z^n q^{n(n-1)/2}
=\prod_{n=1}^\infty 
(1-zq^{n-1})(1-q^n/z)(1-q^n),
\label{eq:tripleprod}
\end{equation}
where $z\ne0$ and $|q|<1$.
The Triple Product Identity is originally due to 
Jacobi \cite[Vol I]{jacobi:works}.
Andrews \cite{andrews84} and Lewis \cite{lewis84} have found 
nice combinatorial
proofs. The triple product occurs frequently in the theory of partitions.
For instance, most proofs of the Rogers-Ramanujan identity
crucially depend on the triple product identity. 

{\tt tripleprod(z,q,T)} returns the $q$-series expansion to order $O(q^T)$ of 
Jacobi's triple product \eqn{tripleprod}.
This expansion is found by simply truncating the right side of
\eqn{tripleprod}.

\sn
{\tt
$>$ \quad  tripleprod(z,q,10);
$$
{\frac {{q}^{21}}{{z}^{6}}}-{\frac {{q}^{15}}{{z}^{5}}}+{
\frac {{q}^{10}}{{z}^{4}}}-{\frac {{q}^{6}}{{z}^{3}}}+{
\frac {{q}^{3}}{{z}^{2}}}-{\frac {q}{z}}+1-z+{z}^{2}q-{z}^{
3}{q}^{3}+{z}^{4}{q}^{6}-{z}^{5}{q}^{10}+{z}^{6}{q}^{15}
$$
$>$ \quad  tripleprod(q,q{\carrot}3,10);   
$$
{q}^{57}-{q}^{40}+{q}^{26}-{q}^{15}+{q}^{7}-{q}^{2}+1-q+{q}
^{5}-{q}^{12}+{q}^{22}-{q}^{35}+{q}^{51}
$$
}
The last calculation is an illustration of Euler's Pentagonal
Number Theorem \cite[Cor. 1.7 p.11]{andrews:book}:
\begin{equation}
\prod_{n=1}^\infty (1-q^n)
=\prod_{n=1}^\infty (1-q^{3n-1})(1-q^{3n-2})(1-q^{3n})
=\sum_{n=-\infty}^\infty (-1)^n q^{n(3n-1)/2}.
\label{eq:euler}
\end{equation}

\subsection{The Quintuple Product Identity}
\label{sec-quinprod}

The following identity is the Quintuple Product Identity:
\begin{align}
& (-z,q)_{{\infty }} (-{\frac {q}{z}},q)_{{\infty }} ({z}^{2}q,{q}^{2})
_{{\infty }} ({\frac {q}{{z}^{2}}},{q}^{2})_{{\infty }} (q,q)_{{
\infty }} \label{eq:quinprod}\\
&\quad =\sum _{m=-\infty }^{\infty }\left (\left (-z\right )^{-3\,m}
-\left (-z\right )^{3\,m-1}\right )
{q}^{{\frac {m\left (3\,m+1\right )
}{2}}}.
\nonumber             
\end{align}
Here $|1q|<1$ and $z\ne0$. This identity is the $BC_1$ case of the
Macdonald identities \cite{macdonald72}. The quintuple product identity was
is usually attributed to Watson \cite{watson29}. However it can
be found in Ramanujan's lost notebook
\cite[p. 207]{ramlost}. Also see \cite{berndt93} for more history
and some proofs.

The function \t{quinprod(z,q,T)} returns the quintuple product
identity in different forms:
\begin{enumerate}
\item[(i)] If $T$ is a positive integer it returns the
$q$-expansion of the right side of \eqn{quinprod} to order $O(q^T)$.
\item[(ii)] If $T=$ {\it prodid} then 
 {\tt quinprod(z,q,prodid)} returns the quintuple product identity 
 in product form.
 \item[(iii)] If $T =$ {\it seriesid} then 
 {\tt quinprod(z,q,seriesid)} returns the quintuple product identity 
 in series form.
\end{enumerate}   
\sn
{\tt
$>$ \quad  quinprod(z,q,prodid);
\begin{eqnarray*}
\lefteqn{{{\rm \ }( - z, \,q)_{\infty }}\,{{\rm \ }( - 
{\displaystyle \frac {q}{z}} , \,q)_{\infty }}\,{{\rm \ }(z^{2}\,q, \,q^{2})_{\infty }}\,{{\rm \ }({\displaystyle \frac {q}{z^{2}}} , 
\,q^{2})_{\infty }}\,{{\rm \ }(q, \,q)_{\infty }}=} \\
 & & {({\displaystyle \frac {q^{2}}{z^{3}}} , \,q^{3})_{
\infty }}\,{{\rm \ }(q\,z^{3}, \,q^{3})_{\infty }}\,{{\rm \ }(q^{
3}, \,q^{3})_{\infty }} + z\,{{\rm \ }({\displaystyle \frac {q}{z^{3}}
} , \,q^{3})_{\infty }}\,{{\rm \ }(q^{2}\,z^{3}, \,q^{3})_{\infty }}\,{{\rm \ }(q^{3}, \,q^{3})_{\infty }}
\end{eqnarray*}
$>$ \quad  quinprod(z,q,seriesid);
\begin{eqnarray*}
\lefteqn{{{\rm \ }( - z, \,q)_{\infty }}\,{{\rm \ }( - 
{\displaystyle \frac {q}{z}} , \,q)_{\infty }}\,{{\rm \ }(z^{2}\,q, \,q^{2})_{\infty }}\,{{\rm \ }({\displaystyle \frac {q}{z^{2}}} , 
\,q^{2})_{\infty }}\,{{\rm \ }(q, \,q)_{\infty }}=} \\
 & & {\displaystyle \sum _{m= - \infty }^{\infty }} \,(( - z)^{(
 - 3\,m)} - ( - z)^{(3\,m - 1)})\,q^{(1/2\,m\,(3\,m + 1))}
\mbox{\hspace{8pt}}
\end{eqnarray*}
$>$ \quad  quinprod(z,q,3);
\begin{eqnarray*}
\lefteqn{(z^{12} + {\displaystyle \frac {1}{z^{11}}} )\,q^{22} + 
( - z^{9} - {\displaystyle \frac {1}{z^{8}}} )\,q^{12} + (z^{6}
 + {\displaystyle \frac {1}{z^{5}}} )\,q^{5} + ( - z^{3} - 
{\displaystyle \frac {1}{z^{2}}} )\,q + 1 + z} \\
 & & \mbox{} + ( - {\displaystyle \frac {1}{z^{3}}}  - z^{4})\,q
^{2} + ({\displaystyle \frac {1}{z^{6}}}  + z^{7})\,q^{7} + ( - 
{\displaystyle \frac {1}{z^{9}}}  - z^{10})\,q^{15} + (
{\displaystyle \frac {1}{z^{12}}}  + z^{13})\,q^{26}
\mbox{\hspace{15pt}}
\end{eqnarray*}
}
Let's examine a more interesting application.
Euler's infinite product may be dissected according to the residue
of the exponent of $q$ mod $5$:
$$
\prod_{n=1}^\infty (1-q^n)
= E_0(q) + qE_1(q^5) + q^2E_2(q^5)
+ q^3E_3(q^5) + q^4E_2(q^5).
$$
By \eqn{euler} we see that $E_3=E_4=0$ since 
$n(3n-1)/2\equiv0$, $1$ or $2$ mod $5$. Let's see if we can identify
$E_0$.
\sn{\tt
$>$ \enskip with(qseries): \ssn
$>$ \enskip EULER:=etaq(q,1,500): \ssn
$>$ \enskip E0:=sift(EULER,q,5,0,499);

\begin{eqnarray*}
\lefteqn{{\it E0} := 1 + {q} - {q}^{3} - {q}^{7} - {q}^{8} - {q}
^{14} + {q}^{20} + {q}^{29} + {q}^{31} + {q}^{42} - {q}^{52} - {q
}^{66}} \\
 & & \mbox{} - {q}^{69} - {q}^{85}  + {q}^{99}\mbox{\hspace{284pt}}
\end{eqnarray*}

$>$ \enskip jacprodmake(E0,q,50);
\[
{\displaystyle \frac {{\rm JAC}(\,2, 5, { \infty}\,)\,{\rm JAC}(
\,0, 5, { \infty}\,)}{{\rm JAC}(\,1, 5, { \infty}\,)}}
\]
$>$ \enskip jac2prod(");
\[
{\displaystyle \frac {{{\rm \:}(\,{q}^{5}, {q}^{5}\,)_{{ \infty}}
}\,{{\rm \:}(\,{q}^{2}, {q}^{5}\,)_{{ \infty}}}\,{{\rm \:}(\,{q}
^{3}, {q}^{5}\,)_{{ \infty}}}}{{{\rm \:}(\,{q}, {q}^{5}\,)_{{ 
\infty}}}\,{{\rm \:}(\,{q}^{4}, {q}^{5}\,)_{{ \infty}}}}}
\]
$>$ \enskip quinprod(q,q\carrot5,20):\ssn
$>$ \enskip series(",q,100);
$$
1+q-{q}^{3}-{q}^{7}-{q}^{8}-{q}^{14}+{q}^{20}+{q}^{29}+{q}
^{31}+{q}^{42}-{q}^{52}-{q}^{66}-{q}^{69}-{q}^{85}+O\left (
{q}^{99}\right )
$$
}
From our \maple session it appears that
\begin{equation}
E_0 =  \frac{ (q^5;q^5)_\infty (q^2;q^5)_\infty (q^3;q^5)_\infty}
{ (q;q^5)_\infty (q^4;q^5)_\infty},
\label{eq:e0id}
\end{equation}
and that this product can be gotten by replacing $q$ by $q^5$
and $z$ by $q$ in the product side of the quintuple product
identity \eqn{quinprod}.

\EX{12}{(i) \quad Use the quintuple product identity \eqn{quinprod}
and Euler's pentagonal number theorem to prove \eqn{e0id} above.\ssn
(ii) \quad Use \maple to identify and prove product expressions for $E_1$
and $E_2$. \ssn
(iii) \quad This time see if you can repeat (i), (ii) but
split the exponent mod $7$.\ssn
(iv) \quad  Can you generalize these results to arbitrary modulus? 
}

\subsection{Winquist's Identity}
\label{sec-winquist}

Back in 1969, Lasse Winquist \cite{winquist} discovered a remarkable
identity
\begin{align}
\label{eq:winquist}
& (a;q)_\infty  (q/a;q)_\infty  (b;q)_\infty  (q/b;q)_\infty (ab;q)_\infty
(q/(ab);q)_\infty  (a/b;q)_\infty  (b/(aq);q)_\infty  (q;q)_\infty^2
 \\
\nonumber
&\quad
=
\sum_{n=0}^\infty\sum_{m=-\infty}^\infty
            (-1)^{n+j}(
                (a^{-3n}-a^{3n+3})(b^{-3m}-b^{3m+1}) \\
&\qquad\qquad +
                (a^{-3m+1}-a^{3m+2})(b^{3n+2}-b^{-3n-1}))
                q^{3n(n+1)/2+m(3m+1)/2}.
\nonumber                 
\end{align}
By dividing both sides by $(1-a)(1-b)$ and letting $a$, $b\rightarrow1$
he was able to express the product $\displaystyle\prod_{n=1}^\infty
(1-q^n)^{10}$ as a double series and prove Ramanujan's
partition congruence
$$
p(11n+6)\equiv0\pmod{11}.
$$
This was probably the first truly elementary proof of Ramanujan's
congruence modulo $11$. The interested reader should see Dyson's
article \cite{dyson72} for some fascinating history
on identities for powers of the Dedekind eta function and how 
they led to the Macdonald identities. A new proof of Winquist's
identity has been found recently by S-Y Kang \cite{kang}.
Mike Hirschhorn \cite{hirschhorn87} has found a four-parameter
generalization of Winquist's identity.

The function \t{winquist(a,b,q,T)} returns the 
$q$-expansion of the right side of \eqn{winquist} to order $O(q^T)$.

We close with an example.
For $1< k < 33$ define
$$
Q(k)= \prod_{n=1}^\infty (1-q^{k})(1-q^{33-k})(1-q^{33}).
$$
Now define the following functions:
$$
A_0=Q(15),\quad
A_3=Q(12),\quad
A_7=Q(6),\quad
A_8=Q(3),\quad
A_9=Q(9);
$$
\begin{align*}
B_0&=Q(16) - q^2 Q(5),\\
B_1&=Q(14) - q Q(8),\\
B_2&=Q(13) - q^3 Q(2),\\
B_4&=Q(7) + q Q(4),\\
B_7&=Q(10) + q^3 Q(1)  .
\end{align*}
These functions occur in Theorem 6.7 of \cite{garvan88} as well
as the function $A_0B_2 - q^2 A_9B_4$.
\sn{\tt
$>$ \enskip with(qseries): \ssn
$>$ \enskip Q:=n->tripleprod(q\carrot{n},q\carrot33,10):\ssn
$>$ \enskip A0:=Q(15): \enskip   A3:=Q(12): \enskip  A7:=Q(6): \enskip   A8:=Q(3): \enskip  A9:=Q(9): \ssn
$>$ \enskip B2:=Q(13)-q\carrot3*Q(2): \enskip B4:=Q(7)+q*Q(4): \ssn
$>$ \enskip IDG:=series( ( A0*B2-q\carrot2*A9*B4),q,200): \ssn
$>$ \enskip series(IDG,q,10);
$$
1-{q}^{2}-2\,{q}^{3}+{q}^{5}+{q}^{7}+{q}^{9}+O\left ({q}^{
11}\right )
$$
$>$ \enskip jacprodmake(IDG,q,50);
\[
{\displaystyle \frac {{\rm JAC}(\,2, 11, { \infty}\,)\,{\rm JAC}(
\,3, 11, { \infty}\,)^{2}\,{\rm JAC}(\,5, 11, { \infty}\,)}{{\rm 
JAC}(\,0, 11, { \infty}\,)^{3}}}
\]
$>$ \enskip jac2prod(");
\begin{eqnarray*}
\lefteqn{{{\rm \:}(\,{q}^{2}, {q}^{11}\,)_{{ \infty}}}\,{{\rm \:}
(\,{q}^{9}, {q}^{11}\,)_{{ \infty}}}\,{{\rm \:}(\,{q}^{11}, {q}^{
11}\,)_{{ \infty}}}\,{{\rm \:}(\,{q}^{3}, {q}^{11}\,)_{{ \infty}}
}^{2}\,{{\rm \:}(\,{q}^{8}, {q}^{11}\,)_{{ \infty}}}^{2}\,{{\rm 
\:}(\,{q}^{5}, {q}^{11}\,)_{{ \infty}}}} \\
 & & {{\rm \:}(\,{q}^{6}, {q}^{11}\,)_{{ \infty}}}
\mbox{\hspace{245pt}}
\end{eqnarray*}
$>$ \enskip series(winquist(q\carrot5,q\carrot3,q\carrot11,10),q,20);
\[
1 - {q}^{2} - 2\,{q}^{3} + {q}^{5} + {q}^{7} + {\rm O}(\,{q}^{9}
\,)
\]
$>$ \enskip series(IDG-winquist(q\carrot5,q\carrot3,q\carrot11,10),q,60);
\[
{\rm O}(\,{q}^{49})
\]
}
From our \maple session it seems that
\begin{align}
A_0 B_2 - q^2 A_9 B_ 4 &=
\qinf{2}{11}
\qinf{9}{11}
\qinf{11}{11}
\qinf{3}{11}^2
\qinf{8}{11}^2\label{eq:id11}\\
&\quad \qinf{5}{11}
\qinf{6}{11},
\nonumber      
\end{align}
and that this product appears in Winquist's identity on replacing
$q$ by $q^{11}$ and letting $a=q^5$ and $b=q^3$.

\EX{13}{(i) \quad Prove \eqn{id11} by using the triple product identity
\eqn{tripleprod}
to write the right side of Winquist's identity \eqn{winquist} as
a sum of two products.\ssn
(ii)\quad In a similar manner find and prove a product form for
$$
A_0B_0 - q^3 A_7 B_4.
$$
}

\bibliographystyle{plain}
\bibliography{refs}

\end{document}